\documentclass[11pt]{amsart}

\usepackage[english]{babel}
\usepackage{enumitem}
\usepackage[ansinew]{inputenc}
\usepackage{amsmath}
\usepackage{hyperref}
\usepackage[T1]{fontenc}
\usepackage{lmodern}
\usepackage{amssymb}
\usepackage{amscd}
\usepackage{marginnote}
\usepackage[all]{xy} % XY-Matrix
\usepackage{dsfont}%einheitsmatrix
\usepackage{stmaryrd} % Eckige Klammern [[ \llbracket \rrbracket ]]

\setenumerate[1]{label=(\alph*)}
\setenumerate[2]{label=(\roman*)}

\theoremstyle{plain}
\newtheorem{theorem}{Theorem}[section]
\newtheorem*{theorem*}{Theorem}
\newtheorem{lemma}[theorem]{Lemma}

\newtheorem*{conjecture*}{Conjecture}

\theoremstyle{remark}
\newtheorem{remark}[theorem]{Remark}

\theoremstyle{definition}
\newtheorem{example}[theorem]{Example}

\newcommand{\glzw}[1]{GL_2 (#1)}
\newcommand{\gln}[1]{GL_n (#1)}
\newcommand{\rhob}{\bar{\rho}}
\newcommand{\Dr}{D_{\rhob}}
\newcommand{\Cr}{\mathcal{C}}
\newcommand{\img}{\text{im}}
\newcommand{\Hom}{\text{Hom}}
\newcommand{\Ad}{\text{Ad}}
\newcommand{\Mn}{\text{M}_n}
\newcommand{\rhosign}{\rho_{sign}}
\newcommand{\Fp}{\mathbb{F}_p}
\newcommand{\ZpT}[1]{\Zp \left\llbracket #1  \right\rrbracket}
\newcommand{\Zp}{\mathbb{Z}_p}
\newcommand{\idM}{\mathds{1}}
\newcommand{\act}[2]{{i(#1)}{#2}{i(#1)}^{-1}}
\newcommand{\xym}[1]{\begin{xy}\xymatrix{#1}\end{xy}}
\newcommand{\xymap}[1]{\begin{xy}\xymatrix@R=7pt{#1}\end{xy}}
\newcommand{\maxI}{\mathfrak{m}}
\newcommand{\Rr}{\mathcal{R}}
\newcommand{\Krdim}{\mathrm{Krulldim}}
\newcommand{\runiv}{\mathcal{R}_\text{univ}(\rhob)}
\newcommand{\vphi}{\varphi}

\hypersetup{
pdfauthor = {J. Sprang},
pdftitle = {A universal deformation ring with unexpected Krull dimension}
}

\begin{document}

\title{A universal deformation ring with unexpected Krull dimension}

\author{Johannes Sprang}

%\adress{Johannes Sprang \at  Universit\"at Regensburg, 93040 Regensburg, Germany \\ \email{johannes.sprang@mathematik.uni-regensburg.de }}

\maketitle
%\subjclass[2010]{11F80}
\begin{abstract}
	A well known result of B. Mazur gives a lower bound for the Krull dimension of the universal deformation ring associated to an absolutely irreducible residual representation in terms of the group cohomology of the adjoint representation. The question about equality - at least in the Galois case - also goes back to B. Mazur. In the general case the question about equality is the subject of Gouv\^{e}a's ``Dimension conjecture''.
	In this note we provide a counterexample to this conjecture. More precisely, we construct an absolutely irreducible residual representation with smooth universal deformation ring of strict greater Krull dimension as expected.
	%In this note we provide a counterexample to a conjecture due to F. Gouv\^{e}a, which says that the Krull dimension of the universal deformation ring - as defined by B. Mazur - associated to an absolutely irreducible residual representation could be expressed by the dimension of the group cohomology of the adjoint representation.

\end{abstract}

\section[Introduction]{Introduction}

Motivated by H. Hida's work \cite{hida1} on certain one-parameter families of $p$-adic Galois representations attached to ordinary $p$-adic cusps forms, B. Mazur lay in his seminal paper \cite{mazur1} the cornerstone for a systematic study of such families by founding the deformation theory of residual representations of pro-finite groups. In many cases all such deformations of a residual representation to complete Noetherian local rings can be parametrized by a universal deformation ring. So if we want to understand ``how many'' lifts exist to a given residual representation, it is very natural to ask about the size of the universal deformation ring. This question is as old as the deformation theory of residual representations itself. Already B. Mazur observed in his fundamental paper \cite{mazur1} that we can give a lower bound for the Krull dimension of the universal deformation ring in terms of the group cohomology of the adjoint representation. The question about equality is the subject of the ``Dimension conjecture'' due to F. Gouv\^{e}a as stated in \cite[Conjecture 1]{gouvea2} and \cite[S. 287]{gouvea1}:

\begin{conjecture*}[Dimension conjecture]
 Let $\rhob: G\rightarrow \gln{k}$ be an absolutely irreducible residual representation of a pro-finite group $G $ such that the universal deformation ring $\runiv$ exists. Then we have
 \[
  \Krdim(\runiv/p\cdot\runiv)=h_1-h_2
 \]
 with
 \[
  h_1:=\dim_{k} H^1(G,\Ad(\rhob))\quad\text{and}\quad h_2:=\dim_{k} H^2(G,\Ad(\rhob)).
 \]
 Here $\Ad(\rhob)$ denotes the associated adjoint representation given by the action of $\rhob$ on $\Mn(k)$ by conjugation.
\end{conjecture*}

In this note we will give a counterexample to this conjecture. The idea is the following: In \cite{zubkov1} A.N. Zubkov shows that for $p\neq 2$ every closed pro-$p$ subgroup of $\glzw{R}$ for an arbitrary pro-finite ring $R$ satisfies some pro-$p$ identity. If we take an appropriate  representation of a pro-finite group $G$ satisfying this relation, there should be hope that this relation is not seen in the universal deformation ring but should be seen by cohomology.\par
By far the most interesting examples in the deformation theory of residual representations are given by representations of Galois groups of number fields with restricted ramification. As already pointed out by B. Mazur in \cite{mazur1}, the Dimension conjecture for Galois representations should be seen as a vast generalization of Leopoldt's conjecture. For some results and evidence for the above conjecture in the case of (tame) Galois representations we refer the interested reader to G. B\"ockle's paper \cite{boeckle1}.\par
Since our counterexample is constructed from an abstract group it does not directly touch the question about the Krull dimension of the deformation ring in the case of Galois representations, but in a certain sense it shows which problems have to be taken care of when dealing with the case of Galois representations.\par

\section*{Acknowledgement}
The idea for this counterexample arose out of my Diploma thesis under the mentoring of Prof. N. Naumann. I would like to use this opportunity to thank Prof. N. Naumann for introducing me to the deformation theory of Galois representations and for all mentoring advice during the development of my thesis.

\section[The Deformation functor]{The Deformation functor}

In this section we recall some basic results from the deformation theory of residual representations. Let $p$ be a prime and $k$ be a finite field of characteristic $p$. Fix a continuous representation $\rhob : G \rightarrow \gln{k}$ of a pro-finite group $G$ satisfying the following finiteness condition:
\begin{itemize}
 \item[($\Phi_p$)] The pro-$p$-completion of every open subgroup $G_0$ of $G$ is topological finitely generated.
\end{itemize}
Such a representation $\rhob$ will be called \emph{residual representation}.
We consider the category $\mathcal{C}$ of complete Noetherian local rings with residue field $k$. The morphisms of $\mathcal{C}$ are homomorphisms of local rings which induce the identity on the residue field $k$. A \emph{lift} of $\rhob$ to $R\in\Cr$ is a continuous homomorphism $\rho:G\rightarrow \gln{R}$ whose reduction modulo the maximal ideal $\maxI_R$ of $R$ coincides with $\rhob$. Two lifts of $\rhob$ are called \emph{strict equivalent} if they are conjugate by an element $\gamma\in \Gamma(R):=\ker(\gln{R}\rightarrow\gln{k})$. A strict equivalence class of lifts is called \emph{deformation}. The \emph{deformation functor} $\Dr(R)$ on $\mathcal{C}$ assigns to every complete Noetherian local ring $R$ the set of all deformations of a given residual representation $\rhob$. The following result is due to B. Mazur:
\begin{theorem}[B. Mazur, {\cite{mazur1}}]
 If $\rhob:G\rightarrow \gln{k}$ is absolutely irreducible and $G$ satisfies the finiteness condition $(\Phi_p)$, then $\Dr$ is representable.
\end{theorem}
For such absolutely irreducible $\rhob$ the representing object $\runiv\in\Cr$ will be called \emph{universal deformation ring}. 

\subsection[Deformation of tame representations]{Deformation of tame representations}\label{chap_tame_def}

If the order of the image of $\rhob$ is prime to $p$ the residual representation is called \emph{tame}. In the tame case N. Boston developed in \cite{boston1} group theoretical techniques, which allow in many cases an explicit calculation of the universal deformation ring. We will recall some of these techniques: In the following we will fix a non-trivial tame residual representation $\rhob : G \rightarrow \gln{k}$ of a pro-finite group $G$ satisfying the finiteness condition ($\Phi_p$). Since $\Gamma(R) $ is a pro-$p$ group for every $R\in\Cr$, every lift $\rho$ of $\rhob$ to $R$ will factor through $G/N$ with  $N$ defined via the maximal pro-$p$ quotient $\ker \rhob/N $ of $\ker \rhob$. For studying all deformations of $\rhob$ we may thus assume in the following without loss of generality that $P:=\ker \rhob$ is a pro-$p$ group. Set $A:=\img \rhob$, then the pro-finite version of Schur-Zassenhaus allows us to choose a lift $\sigma_{W(k)}:A  \hookrightarrow \gln{W(k)}$ of the inclusion $A \hookrightarrow \gln{k}$ of the image of $\rhob$ to the ring of Witt vectors $W(k)$ of $k$. Because $W(k)$ is initial in $\Cr$, the choice of $\sigma_{W(k)}$ yields compatible lifts $\sigma_R:A \hookrightarrow \gln{R}$ for every $R\in \Cr$. Using these lifts we let act $A$ on $\Gamma(R)$ and $P$ by conjugation. The set of all continuous $A$-invariant homomorphisms from $P$ to $\Gamma(R)$ will be denoted by $\Hom_A(P,\Gamma(R))$. Now we have the following result due to N. Boston:
\begin{theorem}[N. Boston, {\cite[(6.1) Proposition]{boston1}}]\label{thm_tame_boston}
If $\rhob: G\rightarrow \gln{k}$ is absolutely irreducible, then there is a natural equivalence of functors:
\[
	\Hom_A(P,\Gamma(\cdot)) \rightarrow \Dr(\cdot)
\]
\end{theorem}

\section[Auxiliary results: Representations and group theory]{Auxiliary results: Representations and group theory}

For the construction of our counterexample the choice of the image of the residual representation will be essential. A faithful representation of the symmetric group $S_3$ will be a good choice:

\begin{lemma}
 Let $k$ be a field of characteristic different form $2,3$. Up to equivalence there are exactly three (absolutely) irreducible representations of $S_3$ over $k$: The trivial $1$-dimensional representation $\mathbf{1}$, the $1$-dimensional sign representation $\rhosign$ and the faithful two dimensional standard representation $\rho_{std}:S_3\hookrightarrow \glzw{k}$. The adjoint representation $\Ad(\rho_{std})$ of $\rho_{std}$ is isomorphic to $\mathbf{1}\oplus \rhosign \oplus \rho_{std}$.
\end{lemma}
\begin{proof}
	By assumption the order of $S_3$ is prime to the characteristic of the field $k$, thus all representations are semi-simple. Further, the prime field of $k$ is already a splitting field for $S_3$. Now, all assertions follow by looking at characters of $S_3$.
\end{proof}

One of the main ingredients of the counterexample is the following result of A. N. Zubkov:
\begin{theorem}[A. N. Zubkov]
 Let $p\neq 2$ be a prime and $R$ be a commutative pro-finite ring. There is no continuous injective goup homomorphism from a non-abelian free pro-$p$ group onto a closed subgroup of $\glzw{R}$.
\end{theorem}
\begin{proof}
 \cite[Theorem 4.1]{zubkov1} and \cite[Theorem 2.1]{zubkov1}.
\end{proof}

\section[The counterexample]{The counterexample}

In this section we will give a counterexample to the following conjecture of F. Gouv\^{e}a:
\begin{conjecture*}[Dimension conjecture]
 Let $\rhob: G\rightarrow \gln{k}$ be an absolutely irreducible residual representation of a pro-finite group $G $ satisfying the finiteness condition $(\Phi_p)$ with universal deformation ring $\runiv$. Then the Krull dimension of $\runiv/p\cdot\runiv$ is given by 
 \[
  \Krdim(\runiv/p\cdot\runiv)=h_1-h_2
 \]
 with
 \[
  h_1:=\dim_{k} H^1(G,\Ad(\rhob))\quad\text{and}\quad h_2:=\dim_{k} H^2(G,\Ad(\rhob)).
 \]
\end{conjecture*}

Let $p\neq 2,3$ be a prime and fix a faithful representation $i:S_3\hookrightarrow \glzw{\Fp}$. Choose compatible lifts $\sigma_R:S_3\rightarrow\glzw{R}$ as in \ref{chap_tame_def} by using Schur-Zassenhaus and let $S_3$ act on $\Gamma(R)$ by conjugation. Let us consider the ring 
\[ 
	\mathcal{R}:=\ZpT{X_1,X_2,X_3,X_4,Y_1,Y_2,Y_3,Y_4}
\] 
and the elements
\[
		\mathbf{X}:=\idM+\begin{pmatrix}
			X_1 & X_2 \\
			X_3 & X_4
		     \end{pmatrix}
		\text{ and }
		\mathbf{Y}:=\idM+\begin{pmatrix}
			Y_1 & Y_2 \\
			Y_3 & Y_4
		     \end{pmatrix}
\]
in $\Gamma(\mathcal{R})$. Let $P$ be the subgroup of $\Gamma(\mathcal{R})$ topologically generated by all $S_3$-conjugates of $\mathbf{X}$ and $\mathbf{Y}$, i.e.:
\[
	\left\{ \act{a}{\mathbf{X}}: a\in S_3 \right\} \cup \left\{ \act{a}{\mathbf{Y}}: a\in S_3 \right\}.
\]
Thus $S_3$ acts continuously on $P$. We define $G$ as the associated semi-direct product $P \rtimes S_3$. Since $P$ is a topologically finitely generated pro-$p$ group, it is obvious that $G$ satisfies the finiteness condition $(\Phi_p)$.  The absolutely irreducible residual representation 
\[
	\xym{ 
		G \ar@{->>}[r] & S_3 \ar@{^{(}->}[r]^-{i} & \glzw{\Fp}
	}
\]
will be called $\rhob$. Now we can state our counterexample:

\begin{example}
	Let $p\neq 2,3$ be a prime and 
	\[ 
		\rhob:G=P \rtimes S_3\rightarrow \glzw{\Fp}
	\] 
	as constructed above. The universal deformation ring of the absolutely irreducible representation $\rhob$ is given by:
	\[
			\Rr=\ZpT{X_1,X_2,X_3,X_4,Y_1,Y_2,Y_3,Y_4}
	\]
	The universal deformation is induced by the inclusion:
	\[
		G \hookrightarrow \glzw{\Rr}
	\]
	In particular $\Krdim \runiv/p\runiv=8$. On the other hand we have:
	\[
		h_1:=\dim_{\Fp} H^1(G,\Ad(\rhob))=8 \text{ and } h_2:=\dim_{\Fp} H^2(G,\Ad(\rhob))\geq 1
	\]
	and thus:
	\[
		\Krdim \runiv/p\runiv >h_1-h_2
	\]
\end{example}
\begin{proof}
	We use N. Boston's description of the deformation functor for tame representations as $\Hom_{S_3}(P,\Gamma(\cdot))$  as given in Theorem \ref{thm_tame_boston}. Let $F$ be the free pro-$p$ group on two generators $x$ and $y$ and consider the homomorphism $\Phi:F\rightarrow P$ given by $x\mapsto \mathbf{X}$ and $y\mapsto \mathbf{Y}$. We first show:\par
	\emph{\textbf{Step 1:} We have the following bijection, functorial in $R\in\Cr$:
	\[
	\xymap{
		\Hom_{S_3}(P,\Gamma(R)) \ar[r] & \Hom_{\text{pro-}p}(F,\Gamma(R))\\
		\vphi \ar@{|->}[r] & \vphi\circ \Phi
	}
	\]}
	Injectivity is obvious since every $S_3$-invariant homomorphism out of $P$ is determined by the images of $X$ and $Y$. Conversely let $\psi\in  \Hom_{\text{pro-}p}(F,\Gamma(R))$ be given with
	\begin{equation*}
		\psi(x)=\idM+	\begin{pmatrix}
								a_1 & a_2 \\
								a_3 & a_4
							\end{pmatrix} \quad\text{and}\quad
		\psi(y)=\idM+	\begin{pmatrix}
								b_1 & b_2 \\
								b_3 & b_4
							\end{pmatrix}\\
	\end{equation*}
	for $a_i$ and $b_i$ in the maximal ideal $\maxI_R$ of $R$. The ring homomorphism $\mathcal{R}\rightarrow R$ given by $X_i\mapsto a_i$, $Y_i\mapsto b_i$ induces an $S_3$-invariant homomorphism $\varphi:\glzw{\mathcal{R}}\rightarrow \glzw{R}$ with $(\varphi|_{P})\circ \Phi=\psi$.\par
	The functor $\Hom_{\text{pro-}p}(F,\Gamma(\cdot))$ is clearly represented by $\mathcal{R}$. The fact that the universal deformation is given as indicated above follows immediately by making the natural equivalences between  $\Hom_{\text{pro-}p}(F,\Gamma(\cdot))$ and $\Dr(\cdot )$ explicit. In particular, we get $\dim_{\Fp}H^1(G,\Ad(\rhob))=8$, since $H^1(G,\Ad(\rhob))$ is isomorphic to the Zariski tangent space of the functor $\Dr(\cdot)$. It remains to show:\par 
	\emph{\textbf{Step 2:}} $H^2(G,\Ad(\rhob))\neq 0$ \par
	Since the order of $S_3$ is prime to $p$, the Hochschild-Serre spectral sequence yields:
	\begin{equation}\label{dimConj_Bsp_eq1}
		H^2(G,\Ad(\rhob))\cong H^2(P,\Ad(\rhob))^{S_3}
	\end{equation}
	Since $P$ acts trivial on $\Ad(\rhob)$, we get an isomorphism of $\Fp[S_3]$-modules:
	\begin{equation}\label{dimConj_Bsp_eq2}
		H^2(P,\Ad(\rhob))\cong H^2(P,\Fp)\otimes \Ad(\rhob)
	\end{equation}
	Combining \eqref{dimConj_Bsp_eq2} and \eqref{dimConj_Bsp_eq1} we obtain:
	\[
		%H^2(G,\Ad(\rhob))\cong (H^2(P,\Fp)\otimes \Ad(\rhob))^{S_3}\cong \Hom_{\Fp[S_3]}(H^2(P,\Fp)^{*},\Ad(\rhob))
		H^2(G,\Ad(\rhob))\cong \Hom_{\Fp[S_3]}(H^2(P,\Fp)^{*},\Ad(\rhob))
	\]
	As every irreducible $\Fp[S_3]$-module occurs in $\Ad(\rhob)$, and since Zubkov's theorem implies that $H^2(P,\Fp)^{*}$ is non-trivial, we can conclude the non-triviality of $H^2(G,\Ad(\rhob))$.
	
\end{proof}

\begin{remark}
We note that $G$ is a quotient of a group $G'$ such that the residual representation
\[
	\xymap{
		G' \ar@{->>}[r] & G \ar[r]^-{\rhob} & \glzw{\Fp}
	}
\]
has the same universal deformation ring but is cohomological unobstructed, i.e. $\dim_{\Fp} H^2(G',\Ad(\rhob))=0$. So for $G'\rightarrow \glzw{\Fp}$ the Krull dimension of the universal deformation ring is indeed given by $h_1-h_2$. The group $G'$ can be constructed as follows: Let $F$ be the free pro-$p$ group on 12 generators $(x_\sigma)_{\sigma\in S_3}$ and $(y_\sigma)_{\sigma\in S_3}$. We let $S_3$ act on $F$ in the following way. Each $\tau \in S_3$ defines an automorphism of $F$ by:
\[
	x_\sigma \mapsto x_{\tau\sigma},y_\sigma \mapsto y_{\tau\sigma}
\]
This action gives $G'$ as the semi-direct product $F\rtimes S_3$. The map $G'\rightarrow G$ is induced by the unique $S_3$-equivariant map $F\rightarrow P$ which maps $x_1$ to $\mathbf{X}$ and $y_1$ to $\mathbf{Y}$.
\end{remark}

\bibliographystyle{amsalpha} 
\bibliography{DimConj_MZ}

\providecommand{\bysame}{\leavevmode\hbox to3em{\hrulefill}\thinspace}
\providecommand{\MR}{\relax\ifhmode\unskip\space\fi MR }
\providecommand{\MRhref}[2]{%
  \href{http://www.ams.org/mathscinet-getitem?mr=#1}{#2}
}
\providecommand{\href}[2]{#2}
\begin{thebibliography}{Maz89}

\bibitem[B{\"o}c98]{boeckle1}
Gebhard B{\"o}ckle, \emph{The generic fiber of the universal deformation space
  associated to a tame {G}alois representation}, Manuscripta Math. \textbf{96}
  (1998), no.~2, 231--246. \MR{1624537 (99j:11057)}

\bibitem[Bos91]{boston1}
Nigel Boston, \emph{Explicit deformation of {G}alois representations}, Invent.
  Math. \textbf{103} (1991), no.~1, 181--196. \MR{1079842 (91j:11041)}

\bibitem[Gou95]{gouvea2}
Fernando~Q. Gouv{\^e}a, \emph{Deforming {G}alois representations: a survey},
  Seminar on {F}ermat's {L}ast {T}heorem ({T}oronto, {ON}, 1993--1994), CMS
  Conf. Proc., vol.~17, Amer. Math. Soc., Providence, RI, 1995, pp.~179--207.
  \MR{1357212 (96j:11071)}

\bibitem[Gou01]{gouvea1}
\bysame, \emph{Deformations of {G}alois representations}, Arithmetic algebraic
  geometry ({P}ark {C}ity, {UT}, 1999), IAS/Park City Math. Ser., vol.~9, Amer.
  Math. Soc., Providence, RI, 2001, Appendix 1 by Mark Dickinson, Appendix 2 by
  Tom Weston and Appendix 3 by Matthew Emerton, pp.~233--406. \MR{1860043
  (2003a:11061)}

\bibitem[Hid86]{hida1}
Haruzo Hida, \emph{Galois representations into
  {${\mathrm{GL}}_2({\mathbf{Z}}_p[[X]])$} attached to ordinary cusp forms},
  Invent. Math. \textbf{85} (1986), no.~3, 545--613. \MR{848685 (87k:11049)}

\bibitem[Maz89]{mazur1}
B.~Mazur, \emph{Deforming {G}alois representations}, Galois groups over {${\bf
  Q}$} ({B}erkeley, {CA}, 1987), Math. Sci. Res. Inst. Publ., vol.~16,
  Springer, New York, 1989, pp.~385--437. \MR{1012172 (90k:11057)}

\bibitem[Zub87]{zubkov1}
A.~N. Zubkov, \emph{Non-abelian free pro-p-groups cannot be represented by
  2-by-2 matrices}, Siberian Mathematical Journal \textbf{28} (1987), 742--747,
  10.1007/BF00969315.

\end{thebibliography}
\end{document}